\documentclass{article}
\usepackage[english]{babel}
\usepackage{amssymb, amsmath}

\newtheorem{scheme}{Scheme}
\newtheorem{prop}{Proposition}

\begin{document}

\title{Simulating elliptic diffusions and orthogonal invariance}

\author{Charles Curry}

\maketitle

\begin{abstract}
We study numerical methods for simulating diffusions that exploit the orthogonal invariance of the Gaussian law, building on the work of Cruzeiro, Malliavin and Thalmeier.
\end{abstract}

\section{Introduction}

Suppose we are given an elliptic operator on $\mathbb{R}^d$
\[ L = \frac{1}{2} a^{i j} (t, x) \frac{\partial^2}{\partial x_i \partial x_j}
   + b^i (t, x) \frac{\partial}{\partial x_i}, \]
and wish to sample from the associated diffusion, which we can interpret as
the probability measure on path space $C ([0, \infty], \mathbb{R }^d)$ solving
the Stroock-Varadhan martingale problem for $L$. Standard methods simulate Ito
processes via numerical solution of stochastic differential equations
\begin{equation}\label{SDE}
X_t = A^i (t, x) \circ d B_t + A^0 (t, x) d t,
\end{equation}
where the driving vector fields $A^i$ are square roots of the diffusion
coefficient $a^{i j}$, i.e. $A A^T = a$. The ellipticity condition assures
these exist, but are only unique up to an orthogonal transformation. The
situation is similar on a Riemannian manifold, but here finding an appropriate
square root is typically problematic, and the usual solution has been to
incorporate the ambiguity by lifting to an SDE on the bundle of orthonormal
frames. This amounts to developing the solution of an SDE on $\mathbb{R}^d$ to
the manifold.

In general, strong simulation to order 1 of an Ito process requires the
numerically expensive simulation of iterated integrals $J_{i j}(0,T) = \int_{0 < s
< t < T} d B_s^{i } d B^j_t$ unless the driving fields commute. Cruzeiro, Malliavin and
Thalmeier \cite{CMT04} introduced a method to exploit the orthogonal invariance to bypass this requirement.
The scheme is a variant of the order 1 Milstein scheme
\begin{eqnarray*}
X_{n+1} - X_n &=& A^0(X_n)h + A^k(X_n) \Delta B^k_t \\
&+& \frac{1}{2}\Big(
(A^k \triangleright A^s)(X_n) + A^i(X_n) \langle [A^s,A^i],A^k \rangle_{X_n}
\Big)
\end{eqnarray*}
where the matrix of iterated integrals is replaced by $\frac{1}{2} (\Delta B^i
\Delta B^j - h \delta_{i j})$, as is possible for commuting vector fields, and
where the usual directional derivatives $A^i \triangleright A^s = A_k^i
\frac{\partial}{\partial x_k} A^j$ are replaced by covariant derivatives for
the Levi-Civita connection of the Riemannian metric defined by inverting the
diffusion matrix.

The scheme is closely related to the Milstein scheme on the trivial
orthonormal frame bundle over $\mathbb{R}^d$. Indeed, the aforementioned
metric is defined such that the elliptic operator is the Riemannian Laplacian
for the metric, and the associated Ito processes are Riemannian Brownian
motions. The basic horizontal vector fields that drive the SDE on the frame
bundle commute up to vertical terms, indeed this is a direct consequence of
the torsion-free property of the connection.

The convergence is of a strong order 1.0, but to a weak solution, i.e. the the relationship between the Brownian motion used to generate the path and the Brownian motion associated to the solution of (\ref{SDE}) is not explicit. Davie has studied methods with this type of convergence, using stochastic couplings, noting the similarity with the methods of Cruzeiro, Malliavin and Thalmeier. Indeed, as he observes, it is not in fact necessary to include the extra terms in the covariant derivative to obtain this order 1.0 using the Milstein scheme.

We focus rather on Castell-Gaines methods \cite{CG}, i.e. those employing ODE solvers on truncations of the stochastic exponential Lie series. These have been considered in the context of coupling methods by Flint and Lyons \cite{FL}. We draw similar conclusions, indeed the link between Castell-Gaines methods on the frame bundle and the base manifold becomes particularly transparent in this context. Moreover, such methods can be generalized to manifolds by employing techniques from geometric integration, see Malham \& Wiese \cite{MW}. We conclude by demonstrating their convergence of order 1.0, and discuss in which situations geometric integration may be applied. The structure is as follows:

\begin{enumerate}
  \item We discuss the geometry of elliptic diffusions, explaining how they can be considered Riemannian Brownian motions for a metric derived from the diffusion tensor.
  \item We review the method of Cruzeiro, Malliavin and Thalmeier and some related schemes that exploit orthogonal invariance by lifting to the orthonormal frame bundle.
  \item We discuss strong convergence to weak solutions, introducing some ideas of Davie.
  \item We propose methods of simulating diffusions on manifolds, based on Castell-Gaines methods and geometric integration, and discuss their convergence in the light of the previous sections.
\end{enumerate}

\section{Geometry of a diffusion}
Consider now a non-degenerate elliptic operator on a Riemannian manifold $M$, expressible in local coordinates as
\[ L = \frac{1}{2} a^{i j} (t, x) \frac{\partial^2}{\partial x_i \partial x_j}
   + b^i (t, x) \frac{\partial}{\partial x_i}, \]
Any such operator is in fact of the form $L=\frac{1}{2}\Delta_M + \tilde{b}$, where $\Delta_M$ is the Riemannian Laplacian for the metric defined through the cometric  $g^{ij}=\sigma^T\sigma$, where $\sigma_{ij}=A^i_j$, see Ikeda \& Watanabe \cite[V.4]{IW}. The associated diffusions are then Brownian motions with drift, which relate to Euclidean Brownian motions through the notion of stochastic development.

\subsection{Development}

A Cartan geometry allows us to develop a curve $x_t$ on a manifold $M$ as a path on the model homogeneous space $G/H$, see Sharpe \cite[V.4.15]{Sharpe}. Of particular interest is the affine development where $G=A(n)$ and $H=GL(n)$, where we may identify the space $G/H\simeq\mathbb{R}^n$ with the tangent space of $M$ at $x_0$. In this case a Cartan connection is an affine connection in the sense of Kobayashi \& Nomizu \cite[III.3]{K&N}, and the development may be computed from integral curves of the corresponding linear parallel transport. Indeed, letting $a_t=(x_t, u_t)$ be the horizontal lift of the curve $x_t$ in the frame bundle $L(M)$, it follows that
$$ u_t(\dot{x}_t) = c_t$$

The inverse procedure of finding a curve $x_t$ which develops to a given path $c_t$ has also been called a development, and by the above result can be considered equivalent to computing projections of integral curves of the basic horizontal vector fields. 

In the particular case that $M$ is an $m$-dimensional submanifold of $\mathbb{R}^n$, the induced Riemannian development has an elegant interpretation. Any path $q(t)$ taking values in $\mathbb{R}^m$ can be \emph{developed} onto $M$ - intuitively the procedure consists of choosing an initial point $p_0$ on $M$; the surface  $\mathbb{R}^m\subset\mathbb{R}^n$ is held tangent to $M$ at $p_0$, with $q(0)$ the point of contact. We then roll $\mathbb{R}^m$ along $M$ without slipping or twisting, such that the point of contact between $M$ and $\mathbb{R}^m$ at time $t$ is $q(t)$. 

We obtain by this procedure a curve $p(t)$ that is related to $q(t)$ by $p(t)=g(t)q(t)$, where $g(t)$ is a curve in the Euclidean group $E(n)$. Indeed 
$$g(t)(x) = A(t)\big(x-q(t)\big) + p(t),$$
where $A(t)\in O(n)$. Given $q(t)$, we compute the development $p(t)$ by solving the following coupled system for $A(t),p(t)$:
\[
\left\{
\begin{array}{ll}
\dot{p} = A \dot{q} & \\
\dot{A}A^{-1} u = -B_p(\dot{p},u) & \forall u \in T_{p}(M) \\
\dot{A}A^{-1} v = -B^t_p(\dot{p},v) & \forall v \in T_{p}^{\perp}(M)
\end{array}
\right.
\]
Where $B_x(u,v)$ is the second fundamental form. Note that $\dot{A}A^{-1}$ is the right-trivialized (Darboux) derivative of $A(t)$, and takes values in ${o}(n)$.

\subsection{Riemannian Brownian motion (with drift)}

The notions of development of a path detailed above extend naturally to stochastic developments of a semimartingale, see Emery \cite[8.31]{Emery}. In particular, Brownian motion on a manifold $M$ is characterized as a stochastic process that develops to a Brownian motion in the tangent space. By anologues of the above results, this may be constructed by solving an SDE in the orthonormal frame bundle. 

Indeed, let $A=\frac{1}{2}\Delta_M + b$ be the sum of a Riemannian Laplacian and a first order term $b$. Let $\tilde{L}_i = B(e_i)$ be basic horizontal fields on the orthonormal frame bundle $O(M)$, and $\tilde{L}_0$ the horizontal lift of the vector field $b$. Suppose $r_t$ solves the following SDE on $O(M)$:
\[
dr_t = \tilde{L}_0(r_t)dt + \tilde{L}_i(r_t)dB^i_t,
\]
Then the projection $x_t=\pi(r_t)$ is a Markov process, an $A$-diffusion, and its law depends only on the initial value of $x_0$ (and not $r_0$)

\subsection{Geometry of an elliptic diffusion}
Recall that we introduced a Riemannian metric $g_{ij}$ defined through the vector fields $A_i$. These form a diagonalizing frame; the components of the Levi-Civita connection associated to the Riemannian structure can therefore be computed as
\[
\Gamma^l_{pq} = \frac{1}{2}(K^l_{pq} + K^p_{lq} + K^q_{lp}),
\]
where $K^i_{jk}(x)$ are the structure constants of the Lie algebra generated by the $A_1,\ldots,A_n$, see Olver \cite[12.41]{Olver}. The drift term of the associated diffusion can be read off from the following:
\[
A=\frac{1}{2}\Delta_M + \big(A_0 + g^{ij} \Gamma^k_{ij} \frac{\partial }{\partial x^k}\big)
\]
Note that it is possible to absorb the drift by redefining the connection, at the expense of introducing torsion. We will not pursue this, as the torsion free property of the Levi-Civita connection is important in the sequel. Indeed, 
the horizontal part of the commutator $[B,B']=-2T(B,B')$ of two basic horizontal fields is directly related to the torsion tensor $T$, and hence vanishes in the absence of torsion \cite[III.5.4]{K&N}. This weak commutativity result is the key underlying the CMT schemes. We conclude by giving the local coordinate expression of the basic horizontal fields:

For any $\xi\in\mathbb{R}^n$, the basic horizontal field $B(\xi)$ on a principal bundle with respect to a connection with components $\Gamma^q_{kl}$ in a local coordinate system $(x^i,X^i_j)$ is given by
\[
B(\xi) = X^i_j \xi^j \frac{\partial}{\partial x^i} -
\Gamma^q_{kl} X^l_p X^k_j \xi^j \frac{\partial}{\partial X^q_p},
\] 

\section{CMT schemes}

We restrict ourselves to $\mathbb{R}^n$, so that the frame bundle is naturally trivial, and introduce the standard basis $e^i_j$ of $n\times n$ matrices.
Let $r_t = (x_t,e_t)$ be the solution of the SDE on the orthogonal frame bundle
\[
dr_t = \tilde{L}_0(r_t)dt + \tilde{L}_i(r_t)dB^i_t.
\]
\begin{scheme}
The Milstein scheme on the frame bundle takes the form
$$
r_{t+h} = r_t + \tilde{L}_0(r_t)\Delta t + \tilde{L}_i(r_t)\Delta W^i_t + \big(\tilde{L}_i \triangleright \tilde{L}_j\big)(r_t) J_{ij}(t,t+h)
$$
\end{scheme}
For the particular case of $n=2$, so that $SO(2)$ is one dimensional and can be parametrized by an angle $\theta$, the above scheme has been utilized under the name 3D-$\theta$ scheme, see Abe \& Giles \cite{A&G}.

The projection of the Milstein scheme for $(\hat{x}_t,\hat{e}_t)$ is 
\begin{eqnarray*}
\hat{x}_{t+h} &=& \hat{x}_t + A_0(\hat{x}_t)\Delta t + \hat{e}^l_k(t) A_l(\hat{x}_t)\Delta W^k_t
\\
&&+ \frac{1}{2}\hat{e}(t)^l_k \hat{e}^{l'}_j \big(A_l\circ A_{l'} - \Gamma^i_{l,l'} A_i\big)(\hat{x}_t))\big(\Delta B^k_t \Delta B^j_t - h\delta^k_j\big)
\end{eqnarray*}
The existence of a canonical trivialization of the frame bundle allows us to replace all the $\hat{e}^i_j$ above with a global section of the frame bundle, i.e. $e^i_j$
\begin{scheme}
The CMT scheme takes the form
$$
X_t = X_0 + A_0(X_0)\Delta t + A_i(X_0)\Delta B^i_t + (A_i\blacktriangleright A_j)(X_0)\big(\Delta B^i_t \Delta B^j_t - h\delta^i_j\big)
$$
or equivalently,
\begin{eqnarray*}
X_t  &=& X_0 + A_0(X_0)\Delta t + A_i(X_0)\Delta B^i_t + (A_i\triangleright A_j)(X_0)\big(\Delta B^i_t \Delta B^j_t - h\delta^i_j\big)
\\
&&+ A_i(X_0) K^i_{jk}(X_0)\big(\Delta B^k_t \Delta B^j_t - h\delta^k_j\big)
\end{eqnarray*}
\end{scheme}
The term involving the covariant derivative $\blacktriangleright$ above has been rewritten to use the structure constants directly, this is possible by symmetry.

Note that if the geometry of the diffusion is such that the Levi-Civita connection is flat, then the commutator $[B,B']$ of the basic horizontal fields vanishes, as this is in general directly related to the curvature \cite[III.5.4]{K&N}. In this case, the parallel transport is path-independent, and the $\hat{e}^i_j$ can be replaced by the global parallel transport of the initial frame.
\begin{scheme}
The 2D-$\theta$ method is the Milstein scheme on the flat frame bundle utilizing a global section of the bundle.
\end{scheme}

If we are to consider weak schemes, we require also additional terms from the stochastic Taylor expansion, namely those involving stochastic integrals $J_{0i}(t)$ and $J_{i0}(t)$, defined as  $\int_0^t s dW^i_s$ and $\int_0^t W^i_s ds$ respectively. The following method has weak order 2.0:
\begin{scheme}[Alves-Cruzeiro]
The AC scheme \cite{A&C} takes the form
\begin{eqnarray*}
X_t &=& X_0 + A_0(X_0)\Delta t + A_i(X_0)\Delta B^i_t + (A_i\blacktriangleright A_j)(X_0)\big(\Delta B^i_t \Delta B^j_t - h\delta^i_j\big)
\\
&&+ \frac{1}{2}(L\circ A_0)(X_0) (\Delta t)^2  + (\tilde{L}\circ A_i)(X_0) J_{0i}(t) + (A_i\triangleright A_0)(X_0) J_{i0}(t),
\end{eqnarray*}
where $L$ is the elliptic operator whose associated diffusions we are simulating, and $\tilde{L}$ is the following covariant modification of the operator:
\[
\tilde{L}\circ A_i = \frac{1}{2}A_{\beta}\blacktriangleright (A_{\beta}\blacktriangleright A_i) + A_0\triangleright A_i.
\]

\end{scheme}

\section{Methods weak and strong}
In general, one distinguishes between a weak and a strong solution of stochastic differential equations such as (1) in that one is free to specify a filtered probability space $(\Omega, \mathcal{F}, P)$ and $\mathcal{F}_t$-Brownian motion $B_t$ together with the solution process $X_t$ for weak solutions, whereas a strong solution $X = F(X(0),B)$ must be a measurable functional of the Brownian motion, i.e. essentially generic. The classic example $dX_t = \mathrm{sgn}(X_t)dB_t$ illustrates the distinction potently: for any initial data $X_0$ and Brownian motion $B(t)$, we let $X(t)=X_0 + B_t$ and construct a new Brownian motion $\tilde{B}_t = \int_0^t \mathrm{sgn}(X_s)dB_s$. It can be checked that $X_t$ is a weak solution for $\tilde{B}_t$, but the integral defining $\tilde{B}$ is an Ito integral constructed probabilistically, not pathwise, and we cannot construct $X_t$ from knowledge of the paths of $B_t$ alone. Indeed, no such strong solution of the form $X=F(X_0, B)$ exists for generic initial data. 

This is not to be confused with notions of weak and strong convergence of approximate solutions. Strong convergence typically measures the deviation of simulated sample paths $X^h_t$ from `true' paths $X_t$, whilst weak convergence measures error in functionals $E(f(X_t))$ of the solution. 

\subsection{Strong convergence}
In considering mean square error 
\begin{equation}
E\big((X_t - X^h_t)^2\big)
\end{equation}
 it is often assumed that $X_t$ is a strong solution and that the $X^h=F^h(X_0,B)$ are similarly Brownian functionals. The error is then also a Brownian functional, whose norm we can compute. In practice, this is often estimated by comparing trajectories
 \begin{equation}\label{multilevel}
 E\big((X^{h'}_t - X^h_t)^2\big)
 \end{equation}
 simulated on progressively coarser grids, where the trajectories $X^{h'}$ and $X^h$ are considered approximations of the same underlying Brownian sample $X(\omega)$. This is realized by computing $X(\omega)$ (and possibly its iterated integrals), and using Chen's relations to compute the restrictions to coarser grids. Bounds on (\ref{multilevel}) are important for implementation of variance reduction by multilevel Monte Carlo techniques.

On the other hand, the above set up is more stringent than necessary if we are only interested in accurate simulation of trajectories. In this case, we consider rather the generation of $X_h$ as sampling from a measure $\rho^h$, and compare the the distance from the true diffusion measure $\rho$ using a metric on spaces of measure, typically a Vaserstein $2$-metric. It is precisely this notion of strong error that the original article of Cruzeiro, Malliavin and Thalmeier considers, and that which we refer to by strong convergence in the sequel, unless specified otherwise.

\subsection{Coupling methods}
Davie \cite{Davie} introduced methods similar in spirit to those of Cruzeiro, Malliavin and Thalmeier which exploit couplings between the driving noise terms $B(\omega)$ and alternative noise $B^h(\omega)$. This has been interpreted and extended using the language of rough paths by Flint and Lyons \cite{FL}. Indeed, the $B$ and $B^h$ are lifted to rough paths, and methods of coupling are related to approximations of the rough differential equation driven by the lift of $B^h$.

\subsection{Strong convergence, weak solution}
The key to understanding the CMT methods is that the trajectories $X^h$ we generate are pathwise approximations to weak solutions $(X_t, \tilde{B})$ for a Brownian motion other than the $B_t$ whose increments we generate (it is perhaps better to follow \cite{FL} and say whose rough path signatures we compute, so as to include values of the iterated integrals, even if these are only abelian rough path approximations). 

The relation between $B$ and $\tilde{B}$ above is an equation of the form $\tilde{B}_t= A_t B_t$, where $A_t$ is frame part of the solution on the frame bundle. This can be interpreted as a coupling, but it is not necessary to make the link as long as we are only interested in weak solutions, as if often the case for the study of diffusions.

On the other hand, it is worth noting that in certain cases, coupling methods can give the type of strong convergence necessary for multilevel Monte Carlo, see \cite{Davie}. In future work we aim to study possible results in this direction for CMT schemes, using a rough path approach.

\section{ODE-based methods}
The Milstein-based schemes above rely on the linear structure of $\mathbb{R}^n$ and hence do not readily generalize to manifolds. 

A better strategy is to use an ODE-based method, which separate the stochastic and geometric stages of the simulation by relating trajectories driven by a given Brownian path $B_t$ to integral curves of the exponential Lie series
\[
\Psi_t (x) = \tilde{L}_0(x) \Delta t + \tilde{L}_i(x)\Delta W^i_t + \frac{1}{2}[\tilde{L}_i, \tilde{A}_j](x) (J_{ij}(t) - J_{ji}(t)) + \ldots
\]
Each time step advances a unit of time along the integral curve starting at $X_t$ of a truncation of the above series, i.e. we take $X_{t+h}=y(1)$, where $y(t)$ solves
\[
\dot{y}(t) = \Psi_t (y)
\] 
Suppose that we employ a Lie-Trotter splitting scheme on the above equation, dividing into horizontal and vertical components. As $[\tilde{L}_i, \tilde{L}_j]$ is vertical, the horizontal component depends only on $\Delta W^i_t$, not $J_{ij}$, i.e. we can write
\[
\Psi_t = \Psi^h_t + \Psi^v_t = V^h(\Delta W^i_t) + V^v(\Delta W^i_t, J_{ij}(t) - J_{ji}(t))
\]
The splitting scheme advances by
\begin{align*}
\tilde{y}_{\tau+k} &= \tilde{\Psi}^h_{\tau, \tau+k} \circ \tilde{\Psi}^v_{\tau, \tau+k} \circ \tilde{y}_{\tau} \\
&= \tilde{\Psi}^h_{\tau, \tau+k}(A\Delta W) \circ \tilde{y}_\tau
\end{align*}
where $A$ is the rotation induced in the fibre by following $\Psi^v$. As the splitting scheme is geometric integrator of order 1, the base projection of this order 1.0 Castell Gaines scheme on the frame bundle is equal in law to the order 0.5 scheme.  We have proved the following result
\begin{prop}
The trajectories of the ``strong order 0.5" Castell-Gaines method converge with strong order 1.0.
\end{prop}
This result is similar to those obtained by Davie using methods of coupling.

\section{Global methods}

In the sequel, we consider methods built from from the first order truncation of the exponential Lie series on the frame bundle
$$ \tilde{\Psi}_t^1 = t \tilde{L}_0(X_0) + \tilde{L}_i(X_0)\Delta W^i_t,$$ 
and in the case of weak methods, those incorporating Lie brackets $[\tilde{L}_0, \tilde{L}_i]$. As before, we can distinguish between methods which compute exclusively projections to the manifold $M$ by covariant differentiation of the $A_i$, and those that evolve on the frame bundle. However, further progress requires specification of the geometric integration method used on the ODE, which in turn depends on the geometry of $M$.

\subsection{The sphere $S_2$}

Consider the special case of developing a curve $q(t)\in\mathbb{R}^2$ onto a sphere $S_2$. We embed $S_2\in\mathbb{R}^3$ as the set of points satisfying $x^2 + y^2 + z^2 = 1$; the second fundamental form is then 
\[
B_x(u,v) = x(u^T v),
\]
with transpose
\[
B^t_x(v,w) = v (x^T w) = v||w||
\]
Let $\Omega(t)=\dot{A}(t)A^{-1}(t)$. The development equations can then be written:
\[
\left\{
\begin{array}{ll}
\dot{p} = A \dot{q} & \\
\Omega u = -p\langle \dot{p},u\rangle & \forall\, u\perp p \\
\Omega v = -\dot{p}||v|| & \forall\, v \;||\; p
\end{array}
\right.
\]
\emph{where in addition $A e_z = p$.}\\
\newline
To simplify these equations, we note that left-multiplying by $\Omega$ is equivalent to taking a cross product, i.e. $\exists\omega: \Omega u = \omega\times u$. The equations for $\Omega$ can be rewritten in terms of $\omega$:
\begin{align*}
\omega\times u &= -p\langle \dot{p},u \rangle \\
\omega\times p &= \dot{p}
\end{align*}
The above equations guarantee that $\omega \perp p,\dot{p}$; it is then a simple exercise to check that the unique solution is given by
\[
\omega = \dot{p}\times p.
\]
The map $\mathbb{R}^3\rightarrow so(3)$ that reconstructs $\Omega$ from $\omega$ is traditionally called the hat map and denoted $\Omega = \hat{\omega}$; we have now reduced the development equations to
\begin{align*}
\dot{p} &= A\dot{q} \\
\dot{A} &= \widehat{ (\dot{p}\times p) } A
\end{align*}
We can combine the equations by writing
\[
\dot{p}\times p = A\dot{q}\times Ae_z = \mathrm{sgn}(A) A(\dot{q}\times e_z),
\]
in the oriented case $A\in SO(3)$ we recover the single equation
\[
\dot{A} = \widehat{ A(\dot{q}\times e_z) } A
\]
This reflects the well-known fact that the orthonormal (oriented) frame bundle of $S^2$ is isomorphic to $SO(3)$. Let now (formally) $\dot{q} = (dW^1, dW^2, 0)$. Riemannian Brownian motions on $S^2$ are constructed by $x_t = A_t e_z$, where $A_t$ solves the frame bundle SDE
\begin{equation}\label{S2 SDE}
dA_t = \big( A_2 dW^1_t - A_1 dW^2_t \big)^ {\widehat{}} A_t,
\end{equation}
and $A_i$ are the $i$th columns of $A$. Note that $[\widehat{A}_1, \widehat{A}_2]=\widehat{A}_3$, showing explicitly the weak commutativity property.

\subsection{Diffusions on $S^2$}
To summarize, the proposed method is as follows. 

\begin{itemize}
 \item At each time step, we simulate $(\Delta W^1, \Delta W^2)$.
 \item We solve the Lie group equation
 \[
 \dot{B} = \big( B\cdot[\Delta W^1, \Delta W^2,0]^T \big)^{\widehat{}} B_t,\quad B_0 = A_n
 \] 
and let $A_{n+1}=B_1$.
\item The Lie group equation is solved by a Lie group integrator, the simplest of which is the Lie-Euler method $B_{t+h} = \exp \Big( hB\cdot[\Delta W^1, \Delta W^2,0]^T \big)^{\widehat{}} \Big) B_t$.
\item The trajectories are the projections onto the base, $y_{nk} = A_n e_z$.
\end{itemize}

\subsection{Lie group integration and diffusions}
To employ a Lie group integrator, we need to be able to evolve in the frame bundle using a Lie group action. The case of $S^2$ is particularly nice as the frame bundle is isomorphic to $SO(3)$. Otherwise, we must work on a case by case basis. 

More generally, the development equations for a submanifold of $\mathbb(R)^n$ may be solve using actions of $SE(n)$; on the other hand this will not preserve the evolution on the manifold unless some projections are involved. There is no problem for the case of $S^2$ as knowing the orientation of the frame is enough to know the position on the manifold. In other words, the Gauss map \cite{K&N} is an isomorphism.

It would be interesting to investigate further examples of Castell-Gaines methods on the frame bundle, either by Lie group integration or by projection. This is the topic of further work.

\section{Hypoellipticity}

Suppose the operator 
\[ L = \frac{1}{2} a^{i j} (t, x) \frac{\partial^2}{\partial x_i \partial x_j}
   + b^i (t, x) \frac{\partial}{\partial x_i}, \]
is not elliptic. We may construct a cometric $g^{ij}$ from $a^{ij}$ as before, but this will no longer be invertible to a metric $g_{ij}$. We are in other words in the world of sub-Riemannian geometry. The diffusion is hypoelliptic if a H\"{o}rmander condition holds on $a^{ij}$, i.e. the associated vector fields are bracket generating. By analogy with the work of Davie, we might hope to extend CMT schemes to this context.

For appropriately nice cometrics, there are methods of constructing sub-Riemannian Brownian motions using frame bundles \cite{subRie}, however we will lack control over the torsion and hence be unable to retain the weak commutativity property underlying the CMT schemes. It remains a topic for further work to circumvent this obstacle.

\end{document}